\documentclass[12pt,titlepage,a4paper,draft]{article}

\usepackage{amssymb,latexsym,amsfonts,amsmath,amsthm}

\usepackage{verbatim}

{\makeatletter \@addtoreset{equation}{section}}

\newcommand{\diag}{\mathop{\mathrm{diag}}\nolimits}
\newcommand{\Ker}{\mathop{\mathrm{Ker}}\nolimits}

\newcommand{\R}{{\mathbb{R}}}

\newcommand{\Z}{{\mathbb{Z}}}

\newcommand{\caU}{{\mathcal{U}}}
\newcommand{\mcU}{{\mathcal{U}}}

\newcommand{\h}{ \mathop{ \mathrm{h} {} }\nolimits }
\newcommand{\e}{ \mathop{ \mathrm{e} {} }\nolimits }

\newcommand{\beq}{\begin{equation}}
\newcommand{\eeq}{\end{equation}}

\newcommand{\LR}{\Leftrightarrow}
\newcommand{\Ra}{\Rightarrow}

\newcommand{\bpm}{\begin{pmatrix}}
\newcommand{\epm}{\end{pmatrix}}

\newcommand{\sq}{\hfill \square}

\newtheorem{theorem}{Theorem}[section]
\newtheorem{lemma}[theorem]{Lemma}

\newtheorem{example}[theorem]{Example}

\begin{document}

\title{Hyperinvariant, characteristic and marked   subspaces}

\author{Pudji Astuti
\\
Faculty of Mathematics\\
 and Natural Sciences\\
Institut Teknologi Bandung\\
Bandung 40132\\
Indonesia
         \and
Harald~K. Wimmer\\
Mathematisches Institut\\
Universit\"at W\"urzburg\\
D-97074 W\"urzburg\\
Germany}

\date{\today}

\maketitle

\begin{abstract}

\vspace{2cm}
\noindent
{\bf Mathematical Subject Classifications (2000):}
15A18, 
47A15 
15A57. 
13C99 

\vspace{.2cm}

\noindent
 {\bf Keywords:}  hyperinvariant subspaces,  marked subspaces,
  characteristic sub\-spaces,   invariant subspaces,
Jordan basis.

\vspace{.5cm}
\noindent
{\bf Running title:}  Hyperinvariant  subspaces

\vspace{5cm}
\noindent
{\bf Abstract:} \, 
Let  $ V $ be a finite dimensional vector space over a field 
$K$ and $ f $  a $K$-endomorphism of $V$.
In this paper  we study three types of $f$-invariant subspaces,
 namely hyperinvariant  subspaces, which 
are  invariant under
all endomorphisms of $V$ that commute with $f$, characteristic
subspaces, which  remain fixed under all 
automorphisms of $V$ that commute with $f$, and 
marked subspaces, which have 
a Jordan basis (with respect to 
$f_{|X}$) that can be extended to a Jordan basis of $V$.
 We show that a subspace is  hyperinvariant
if and only if it is characteristic and marked. 
If $ K $ has more than two elements then each 
characteristic subspace is  hyperinvariant.

\end{abstract}

\section{Introduction}

Let
 $ V $ be an $n$-dimensional vector space over a field    $K $
 and
let $f : V \to V $ be   $ K $-linear.
We assume that      the characteristic  
  polynomial
of $f $ splits over $K$ such that all eigenvalues of
$f$ are in $ K$.
In this paper we deal with three types of $ f $-invariant subspaces,
namely  with hyperinvariant, characteristic and 
marked   subspaces.
To describe these three  concepts we use the following
notation. 
Let $ {\rm{Inv}}( V) $ be the lattice of $f$-invariant subspaces of $V$
  and let $  {\rm{End}}_f (V) $  be the algebra of 
all endomorphisms of $V$ that commute with $f$. 
If a  subspace $X $   remains
invariant for all $ g \in  {\rm{End}}_f (V) $
then $X$ is  called
{\em{
 hyperinvariant}} for $f$ \cite[p.\ 305]{GLR}.
Let $  {\rm{Hinv (V) }} $ be the   set of  hyperinvariant 
 subspaces of  $V$. It is obvious that
  $  {\rm{Hinv (V) }} $ 
is a lattice. Because of $ f \in  {\rm{End}}_f (V) $ we have 
$   {\rm{Hinv (V) }} \subseteq {\rm{Inv}}( V ) $.
We refer to  \cite{GLR}, \cite{FHL},   \cite{Lo},  
    \cite{Wu}    
for results on  hyperinvariant   subspaces.
The group of  automorphisms of $V$
that commute with $ f $ will be denoted by 
 ${\rm{Aut}}_f(V) $.
A subspace $ X $ of $V$  
will be   called  
{\em{characteristic}}
(with respect to $f$) 
 if  \,$ X \in  {\rm{Inv}}( V ) $\, and 
\,$ \alpha( X ) = X $\,
for all $  \alpha \in   {\rm{Aut}}_f (V) $.
Let $  {\rm{Chinv (V) }} $ be  set of  characteristic 
 subspaces of  $V$. Obviously, also   $  {\rm{Chinv (V) }} $ 
is a lattice, and \,$   {\rm{Hinv (V) }} 
 \subseteq   {\rm{Chinv (V) }}  $.

Set $ \iota = {\rm{id}}_V $ and $ f^0 =  \iota$.
Let 
\,$ \langle x  \rangle _f  = {\rm{span}} \{ f^i x , i \ge 0 \} $\,
 be  the
cyclic subspace generated by $ x \in V $.
If $ B  \subseteq V$ we define
 \,$
 \langle B \rangle _f =
 \sum _{b \, \in \,  B } \,\,   \langle \,  b  \, \rangle _f $.
Let  $\lambda $ be an eigenvalue of $f$  
such that  $  V_{\lambda} =  \Ker ( f - \lambda  \iota )^n $ 
is the corresponding  generalized eigenspace.  
Let \,$ \dim \Ker  ( f - \lambda \iota )  = k $,
and let  $ s^{t_1} , \dots , s^{t_k} $, be  the elementary divisors
of $ f_{| V_{ \lambda } } $. 
Then there exist vectors $u_1, \dots , u_k$, such that 
\[  V_{\lambda} = 
 \langle u_1 \rangle _{f - \lambda \iota } 
\,  \oplus \, 
\cdots \, 
 \oplus    \langle u_{ k } \rangle _{f - \lambda \iota }, 
\]
and $ ( f -  \lambda    \iota  ) ^{ t_i  -1 } u_i \ne 0 $,
 $ ( f -  \lambda    \iota  ) ^{ t_i } u_i =  0 $, $ i = 1, \dots , k$. 
We call $  U _{ \lambda } =  \{  u_{1}, \dots , u_{k}  \}  $
a set of {\em{generators}} of  $  V_{\lambda} $.
Each  $U _{ \lambda } $ gives rise to a   Jordan basis  of $V _{ \lambda }$,
namely
\begin{multline*}
 \bigl\{ u_{1}, \,  ( f-  \lambda      \iota  )u_{1},
 \, \dots ,   ( f -  \lambda    \iota  ) ^{ t_1 -1 } u_1
 \dots ,
\\
  u_{k}, \,  ( f -  \lambda     \iota  )u_{k} ,
  \, \dots ,   ( f -  \lambda      \iota  ) ^{ t_k -1 } 
 u_{k}
 \bigr\}.
\end{multline*}
Define $ f_{\lambda } = f_{| V_{ \lambda }  }$. Let  $Y $ be an
$ f_{\lambda } $-invariant subspace of  $  V_{\lambda} $.
Then $ Y$ is said to be 
 {\em{marked}}
 in   $  V_{\lambda} $  (with respect to $ f_{\lambda } $)
 if there  exists a set
 $ U_{\lambda} $  
of generators of   $  V_{\lambda} $ 
and  corresponding   integers $ r_{i} $, $  0 \le r_i \le t_i $,
such that
\[ 
 Y  = \,
 \big\langle (f  - \lambda  \iota) ^{ r_{1} } \,  u_{1} \big\rangle
 _{ f -   \lambda  \iota   }  
\, \oplus \, \cdots \, \oplus \, 
  \big\langle (f  - \lambda  \iota) ^{ r_{k} } \,  u_{k} \big\rangle
 _{f -  \lambda \iota     }.  
\]
 Thus $Y$  has a Jordan basis 
which can be extended to a
Jordan basis  of $ V_ {\lambda } $.
Let  $ \sigma (f) = \{\lambda _1,    \dots, \lambda _m  \}$
 be the spectrum of $ f $.
Then
\beq  \label{eq.vc}
V = V_{\lambda_1}  \oplus \dots \oplus  V_{  \lambda_m   }      .
   \eeq 
If  $ X  \in {\rm{Inv}}_f V $
then \,$  X _{  \lambda _i    }  =  X  \cap  V_ { \lambda _i } $\, 
is $ f_{ \lambda _i } $-invariant in  $V_ { \lambda _i  }$,
and 
\beq  \label{eq.xli}
  X = X_{\lambda_1}  \oplus \dots \oplus  X_{  \lambda_m   }      .
\eeq 
We say that  $ X $ is  marked in $V$ if  
 each subspace 
$   X _{  \lambda _i    } $  in \eqref{eq.xli}   
 is  marked in $  V_ { \lambda _i }  $. 
The set of  marked
subspaces of   $V$   will be denoted by ${\rm{Mark (V) }} $.
We assume $  0  \in  {\rm{Mark (V) }} $.
 Marked subspaces can be traced back to \cite[p.\,83]{GLR}.
They
have  been studied in  \cite{Bru}, \cite{FPP}, \cite{AW},
and         \cite{CFP}. For marked $(A,C)$-invariant subspaces
we refer to  \cite{CF} and  \cite{CP}.
We mention applications to
algebraic Riccati equations \cite{AW2} and  to
stability of invariant subspaces of commuting matrices \cite{KP}.

The  following examples 
show  that
to a certain extent
the three types  of invariant subspaces  are   independent of each other.
Suppose $ f $ is nilpotent.
If $ x \in V $ then the smallest nonnegative integer $\ell$
 with
$f^{\ell} x = 0$
is  called
the {\em{exponent}} of $x$. We write  $\e(x) = \ell$.
A  nonzero vector  $x $ 
  is said to have 
\emph{height} $q $ if $x \in f^q V$
and $x \notin f^{q+1} V$.
In this case we write $\h(x) = q$. We set $ \h ( 0 ) = - \infty $.
For \,$j \ge 0 $\, we   define 
\,$  V[ f^j ] =  \Ker f^j $.

\begin{example} \label{ex.fu22}
{\rm{
Let $ K = \Z _2 $.  Consider  $ V = K ^4 $
 and
\[
f = \diag ( 0 ,  N_3) , \:\,
         N_3 = \begin{pmatrix}   0 & 0   & 0
\\
 1 & 0  & 0 \\
                0 & 1 & 0
\end{pmatrix} .
\]
 Let   $e_1, \dots , e_4$,
be the unit vectors of $K^4$. Then  $ f ^3 = 0$ and 
\,$ V = \langle e_1 \rangle _f  \oplus  \langle e_2  \rangle _f  $.
  Define \,$ z =  e_1 + e_3  $\, and
\,$
 Z  =  \langle z     \rangle _f  $.
Then
\[
Z =
\{ 0 ,  z ,  z +e_4 , e_4   \}
=
\big\langle v ; \,\,  \e(v) =  2, \,  \h(v) =  0 ,\, \h(f v) = 2
   \big\rangle _f  .
\]
If
 $ \alpha \in  {\rm{Aut}}_f(V) $ then
 $ |  \alpha (Z ) |   = |    Z  |  $. Moreover $\alpha  $
preserves height and exponent.  Hence
 $ \alpha ( Z ) =  Z $, and
  $ Z $ is characteristic.
 Let
$ g = \diag (1, 0 , 0 ,  0) $   be the orthogonal projection
on $ K e_1 $.
Then $ g  \in  {\rm{End}} _f (V) $.
We have  \,$ g z   = e_ 1 \in g (Z )$, but
$ e_1 \notin Z $.
 Therefore  $ Z $ is not hyperinvariant.
The Jordan bases of $ Z $ are
$J_1 = \{z , e_4 \} $ and
 $J_2 = \{ z  + e_4, e_4 \} $.
If \,$ y \in K^4 $\, then 
 \,$ z  \ne  fy   $\,  and  \, $  z +e_4  \ne fy $.
Hence neither $J_1 $ nor  $J_2 $ can be
extended to a Jordan basis of $ K ^4$. Therefore $ Z $ is not
marked.
}}
\end{example}

\begin{example}
{\rm{
Let
$V = K ^2$  and
\,$
 f = 0 $. 
Then
\,$
  K^2 = \langle  e_1 \rangle _f \oplus \langle  e_2 \rangle _f $ \,
and the subspace
$X = \langle  e_1 \rangle _f  $ is marked.  From
 \,$
\alpha
=  \left( \begin{smallmatrix}      1 &  0 \\ 1 & 1
\end{smallmatrix}    \right)
\in   {\rm{Aut}} _f ( V)
$\,
and $\alpha (e_1 ) = e_1 + e_2 $ follows that
$ X $ is not characteristic.
  }}
\end{example}

In contrast to $  {\rm{Hinv}}(V) $ or  $ {\rm{Chinv}}(V) $ 
the set  $ {\rm{Mark}}(V) $ in general is not a lattice. 

\begin{example} 
{\rm{
$ V = K ^6$, 
$ f = \diag(0,N_3,N_2)$. 
The subspaces  $ Z_1 =  \langle e_5 \rangle $ 
 and $ Z_2 =  \langle e_5 + e_3 + e_1 \rangle $ are marked
but  $ Z_1 + Z_2 =  \langle e_5 \rangle  \oplus 
 \langle  e_3+e_1  \rangle $ is not marked.
Thus the   set of marked subspaces is not closed under addition.
}} 
\end{example}

In this paper we study the following problems.
Under what  conditions  is a marked subspace
 characteristic? When is each characteristic subspace hyperinvariant?
Because of the 
   Lemma \ref{la.ds} below 
 one can
deal separately with  single
  components  $ V_{ \lambda _i } $  in
\eqref{eq.vc} and  the  corresponding  restrictions
$ f _{ \lambda _i} =  f_{| V_{ \lambda _i}} $, $i = 1, \dots , m$.

\begin{lemma} \label{la.ds}
 An $f$-invariant
subspace   $ X  \subseteq V  $ is 
 hyperinvariant (resp.\ characteristic, resp.\  marked) 
if and only if, with respect to $  f _{ \lambda _i} $,
each component $ X _{  \lambda _i } $ in \eqref{eq.xli}
 is   hyperinvariant (resp.\ characteristic, resp.\  marked) 
 in   $   V_{ \lambda _i }$.
\end{lemma}

\noindent Proof. 
If $  \eta  \in   {\rm{End}}_f(V) $ then it is known   (\cite[p.\,223]{Ga})
that the  subspaces $   V_{ \lambda _i }$ in \eqref{eq.vc}
are  invariant under $ \eta $, and
that $\eta _{| V_{ \lambda _i } } \in
   {\rm{End}}_{f _{\lambda _i} }  ( V_{ \lambda _i })$. 
Hence, if  
$  X \in {\rm{Inv}}(   V )  $ then
\eqref{eq.xli}
implies 
\[
  \eta ( X ) =  \eta _{| V_{ \lambda _1 } }
( X _{  \lambda _1 } ) 
\oplus \cdots \oplus 
 \eta _{| V_{ \lambda _i }  }
 ( X _{  \lambda _m } )
\]
Hence if  
\,$X _{  \lambda _i } \in {\rm{Hinv}} (   V_{ \lambda _i } )$,
resp. \,$X _ {  \lambda _i } \in {\rm{Chinv}}(   V_{ \lambda _i } )$,
$ i = 1, \dots ,m $,
then 
 \,$ X \in {\rm{Hinv}}(   V )$, resp.  \,$ X \in {\rm{Chinv}}(   V ).$

 Now suppose now that $ X $ is hyperinvariant.
Let us   show that    $  X_{ \lambda _i } \in 
 {\rm{Hinv}} (   V_{ \lambda _i } ) $, $ i = 1, \dots ,m $.   
Take $ i = 1 $. Set 
$ \hat{V} = 
  V_{ \lambda _2 } \oplus \cdots   \oplus  V_{ \lambda _m }$
and 
$ \hat{X} = 
  X_{ \lambda _2 } \oplus \cdots   \oplus  X_{ \lambda _m }$.
Let 
$ \beta _1 \in 
  {\rm{End}}_{f _{\lambda _1} }  ( V_{ \lambda _1 })  $.
Define  
\,$
  \beta = \beta _1 + {\rm{id}}_{| \hat{V} }  $.
Then $ \beta \in   {\rm{End}}_f  ( V ) $.
Hence 
\,$  \beta (X) \subseteq X =  X_{ \lambda _1 } 
 \oplus   \hat{X} $,
and 
\,$
 \beta (X) =\beta _1 (  X_{ \lambda _1 }  )
\oplus  \hat{X} $.
From    $ X_{ \lambda _1 } \subseteq 
 V_{ \lambda _1  }$ and  
$  \beta _1 (  X_{ \lambda _1 } ) \subseteq 
 V_{ \lambda _1} $ 
we obtain $   \beta _1 (  X_{ \lambda _1 } ) 
\subseteq   X_{ \lambda _1 } $. 
Therefore $  X_{ \lambda _1 }\in   {\rm{Hinv}} (   V_{ \lambda _1 } ) $.
A similar argument shows that 
  \,$X \in {\rm{Chinv}}(   V  )$
implies 
\,$  X_{ \lambda _i } \in {\rm{Chinv}}(   V   _ {\lambda _i }  ) $,
$i = 1, \dots, m$.
In the case of  marked subspaces the assertion is obvious. 
$\sq$

\section{Auxiliary  results}

Because of Lemma \ref{la.ds} it suffices to consider an
endomorphism  $ f $ with only one eigenvalue $\lambda $.
We shall assume
$ \sigma (f ) =  \{ 0 \}$
such that $ f^n = 0 $. 
Let
\beq  \label{eq.elt}
 s^{ t_{1}}, \, 
\dots , \,  s ^{t_{k} } ,  \, \,\,  
 0 < t_{1}  \, \le  \, \cdots  \,   \le  \,  t_{k},
\eeq
be the elementary divisors of $f$.
We call  \,$ U = ( u_1, \dots , u_k) $\,   a
{\em{generator tuple}}
  of $V$ if
\beq \label{eq.fo} 
    V = \langle  u_1  \rangle _f \oplus \cdots \oplus
    \langle  u_k   \rangle _f
\eeq
 and if $U$ is
ordered according to \eqref{eq.elt}  such that
\[
  \e( u_1 )\,  = \, t_1 \,\,   \le \,\,  \cdots \,\, 
 \le  \,\,   \e( u_k ) \,  = \,  t_k.
\]
Let  $ \mathcal{U} $ be the set of   generator tuples  of $V$.
In the following
 we omit  the subscript $f$ in   \eqref{eq.fo}
 and   
we  write \,$ \langle  u_i  \rangle  $\, instead of 
 \,$ \langle  u_i \rangle _f $.
We say that a
$k$-tuple
 $r = (r_1, \dots , r_k) $  of 
integers  is {\em{admissible}}   if
\beq \label{eq.adr}
 0 \leq r_i \leq t_i, \,\, i = 1, \dots , k.
\eeq
Each   $U \in  \mathcal{U}$   together with an
 admissible tuple $r$ gives rise
to a subspace
\beq \label{eq.wru}
W(r, U)
=
 \langle f ^{r_1}  u_1 \rangle  \,  \oplus  \,  \cdots  \,   \oplus
 \,  \langle f ^{r_k}  u_k \rangle  ,
\eeq
which is marked  in $V$. Conversely,
 a subspace $W $ is marked in $V$  
only if $ W = W(r, U) $ for some  $U \in  \mathcal{U}$
   and some admissible  $ r$.
The following example shows that,
in general,
  \,$W(r, U) \neq  W(r, \tilde{U})$\,
if
 $ U \ne \tilde{U} $.
\begin{example}  \label{ex.ex}
{\rm{
Let $V = K^5$ and  $ f = \diag (N_2, N_3) $. 
Then 
\,$ V =  \langle e_1   \rangle   \oplus  \langle e_3 \rangle $,
and $ U = (e_1, e_3) $ and $  \tilde{U} = ( e_1, e_3 + e_1) $ 
are generator tuples. 
 Choose \,$r = (1,0)$.  Then
 the corresponding subspaces
 \,$
  W(r, U) = 
 \langle e_2 \rangle   \oplus  \langle e_3  \rangle $\, 
and \,$  W(r, \tilde{U}) =
\langle e_2  \rangle   \oplus  \langle  e_3 + e_1  \rangle  
 $
 are different from each other.}}
\end{example}

The construction of invariant subspaces of the form
$  W(r, U) $ is a standard procedure in linear algebra
and systems theory. It is used   in
 \cite{Ku}, 
\cite[p.61]{GLR0},   \cite[p.28]{BEGO},
\cite{Ro}.
Hence it is important to know
  whether for a given~$r$   different choices of
$U$ will always  result in the same subspace. 
Theorem~\ref{thm.qv}  will provide a necessary and sufficient condition
for $r$ such that   $  W(r, U) $ is independent of the choice of $U$. 
Let $ r $ be admissible and define  
\beq \label{eq.wr}
   W(r) =
f^{r_1}V \cap V[ f^{t_1 - r_1} ]
 \, + \cdots
+ \,
  f^{r_k}V \cap V[ f^{t_k - r_k} ] .
\eeq
Subspaces  of the form 
 $f^{ \nu  }V $ and 
$  V[f^{\mu  }] $ are  hyperinvariant,
and    \,$   {\rm{Hinv}}(V) $ is a lattice.
Therefore (see e.g. \cite{FHL})  
we have   \,$ W(r ) \in  {\rm{Hinv}}(V)$.

\bigskip

The following lemma shows that
   each $\alpha \in  {\rm{Aut}}_f(V) $ is uniquely determined
by the image of a given  generator tuple.

\begin{lemma} \label{la.bj}
Let    \,$ U = ( u_1, \dots , u_k)
 \in \mathcal{U} $  be given.
For \,$    \alpha \in  {\rm{Aut}}_f(V) $\, define
\,$
\Theta_U( \alpha)    = \bigl( \alpha ( u_1 ) ,   \dots ,  \alpha ( u_k  ) \bigr)$.
{\rm{(i)}}
Then
\[  \alpha \mapsto  \Theta  _U  ( \alpha) , \,\,
 \Theta _U :  {\rm{Aut}}_f(V) \to  \mathcal{U} ,
\]
is a bijection.
{\rm{(ii)}} If
\,$  \tilde{U}  = \Theta (   \alpha ) $\, then
\,$        W(r,   \tilde{U})
= \alpha \bigl(   W(r, U) \bigr)
  $.
\end{lemma}
\noindent Proof.  (i)
It is easy to see that $  \Theta _U( \alpha)    \in  \mathcal{U} $.
Hence $   \Theta _U $ maps  $ {\rm{Aut}}_f(V) $ into  $\mathcal{U}$.
Let  $ x \in V $ and  
\beq \label{eq.ngq} 
x = \sum \nolimits _{i = 1 } ^k  \sum \nolimits _{j = 0 } ^{\e( u_i ) -1 }
c_{i j} f ^j u_i.
\eeq 
Suppose  $ \alpha , \beta \in {\rm{Aut}}_f(V)$ and 
  $  \Theta _U( \alpha)  =   \Theta _U( \beta) = 
(\hat{u} _1  , \dots , \hat{u} _k)  $.
Then 
\[
\alpha (x ) =  \sum \sum c_{i j} f ^j  \hat {u} _i = \beta (x).
\]
Hence $ \alpha  = \beta $, and  $   \Theta _U $ is injective. 
Now consider  
\,$   \tilde{U} = ( \tilde{u}_1, \dots,   \tilde{u} _k  )   \in \mathcal{U}$.
Let  $ x \in V $ be the vector in  \eqref{eq.ngq}.
 Define 
\,$  \gamma : x \mapsto   \sum _i \sum _j \, c_{ij} f ^j  \tilde{u} _i $.
Then $ \gamma \in  {\rm{Aut}}_f(V) $ and $  \tilde{U} =   \Theta _U ( \gamma)$.
Hence   $\Theta _U $ is surjective. 
\\
(ii) It is obvious that
\,$
\alpha \bigl(   W(r, U) \bigr) =
\langle f ^{r_1}   \alpha ( u_1)  \rangle _f \oplus \cdots  \oplus
\langle f ^{r_k}   \alpha  ( u_k ) \rangle _f
 =
  W(r, \tilde{U} )$.
\phantom{00} $\sq$

 In group theory  fully invariant subgroups
play the role of 
 hyperinvariant subspaces.  
Hence the decomposition \eqref{eq.drs}  below 
is an analog to a distributive law in  Lemma~9.3 in \cite[p.\ 47]{FuI}.

\begin{lemma}  \label{la.hdir}
Suppose 
\beq \label{eq.dis}
 V = V_1 \oplus \cdots \oplus V_q, \,\, 
V_i \in  {\rm{Inv}}(V), \, i = 1, \dots , q.
\eeq
{\rm{(i)}}
If $ X $ is a hyperinvariant subspace of $V$, or
\\
{\rm{(ii)}} if  $ X $  characteristic 
and \,$ | K | > 2$, then
\beq  \label{eq.drs}
  X = (  X \cap V_1 ) \oplus \cdots \oplus  (  X \cap V_q ).
\eeq 
\end{lemma} 

\noindent Proof.
If $ x \in V $ then  
\,$
x = \sum \nolimits _{i = 1} ^q x_i, \,\: x _i \in V_i $.
 Set 
$ X _i =  X \cap V_i $,
and  \,$ S =  \oplus _{i = 1 } ^q X_i $. Then  \,$ S \subseteq X$.
To prove the converse inclusion 
we 
note that 
\beq \label{eq.nlf} 
   fx = \sum \nolimits _{i = 1 } ^q  f_{| V_i } (x_i) .
\eeq
(i)  
Let $ \pi _i $ 
be the projection on $V_i$ induced by \eqref{eq.dis}. 
 Then \eqref{eq.nlf} implies  
$\pi _i \in  {\rm{End}}_f(V) $. Hence, if $ x \in X $ then  
and $\pi _i(x)  = x _i \in X $. 
Thus  $ x _i \in X _i  $,
and therefore $ X  \subseteq S $. 
\\
(ii) 
Let $ a \in K $ be different from $0 $ and $1$, and define
$ \gamma _i = \iota - a \pi _i $. Then
 \,$  \gamma _i  \in   {\rm{Aut}}_f(V) $.
Hence 
$  \gamma _i (x)  =  x - a x_i \in X $
if $ x \in X $. Thus we obtain $ x_i \in X _i $.
\phantom{..} 
$ \sq $

\bigskip 

\begin{example} \label{ex.chnd} 
{\rm{In Lemma \ref{la.hdir}(ii)
one can not drop   
the assumption   \mbox{$ | K | > 2 $}.
Suppose  \,$ | K | = 2 $, and  let $ V $ and $f $ be as in
Example~\ref{ex.fu22}. The subspace 
\,$ Z =  \langle e_1 + e_3  \rangle$\, is  
characteristic. 
Both 
 \,$ V_1 = \langle e_1  \rangle$\, and 
\,$ V_2 = \langle e_2  \rangle $\, are in $ {\rm{Inv}}(V) $, 
and we have \,$V  = V_1 \oplus V_2$.
But 
$ Z \cap V_1 = 0$ and $  Z \cap V_2 =   \langle e_4  \rangle$
imply 
$ Z \supsetneqq  ( Z \cap V_1  ) \oplus  (  Z \cap V_2 ) $.
}}
\end{example} 

The next lemma is  an intermediate result.

\begin{lemma} \label{thm.gwa}
Each  hyperinvariant subspace of $V$ is  marked,  and  
\beq \label{eq.inch} 
 {\rm{Hinv}}(V) \subseteq {\rm{Mark}}(V) \cap {\rm{Chin}}(V).
\eeq 
\end{lemma} 

\noindent Proof.   
 Let $ U = (u_1, \dots , u_k)  \in \caU$. If  $X$ is invariant
 then $ X \cap   \langle u_i  \rangle  =  
\langle f^ { r _i }  u_i  \rangle $ for some $r_i$. Thus,  if  $X$
is  hyperinvariant
then  
\eqref{eq.drs} in  Lemma~\ref{la.hdir} implies  
\,$ X = \oplus _{i = 1 }^k   \langle f ^{r _i  } u_i  \rangle$.  
Therefore $X $ is marked,
and \,${\rm{Hinv}}(V) \subseteq  {\rm{Chin}}(V)$\, yields  the 
inclusion \eqref{eq.inch}. 
$ \sq $

\section{Hyperinvariant = characteristic +
marked}

\bigskip 
 We now characterize those marked subspaces
which are characteristic.
 The theorem  below 
 includes results from  \cite{AW2} 
  with
 new proofs.

 \bigskip 
\begin{theorem}  \label{thm.qv}
Let \,$ U \in \mathcal{U} $\, and let $r  = (r_1, \dots , r_k) $\,
be admissible. Then the following statements are equivalent.
\begin{itemize}
\item[\rm{(i)}] The subspace  $  W(r, U) $ is characteristic.
\item[\rm{(ii)}]
The 
subspace
 $  W(r, U) $ is independent of the generator tuple
$U$, i.e.\
\beq \label{eq.fa}   W(r, U)=   W(r, \tilde{U})  \quad  for  \: \:
all \quad    \tilde{U} \in \mathcal{U} .
\eeq
\item[\rm{(iii)}]
The tuples $t = ( t_1, \dots , t_k)$ and  \,$r = (r_1, \dots , r_k) $\,
satisfy
  \beq  \label{eq.r1}
 r_1  \leq  \cdots  \leq  r_k
\eeq
and
\beq  \label{eq.r2}
  t_1 - r_1   \leq   \cdots   \leq  t_k - r_k .
\eeq
\item[\rm{(iv)}]  We have  
 $  W(r, U) = W (r)  $.
\item[\rm{(v)}] 
 $  W (r, U) $ is the unique marked subspace $W$ such that
the elementary divisors   of $W$ and   of \,$V/W$   are
\beq \label{eq.li2}
 s^{t_1 - r_1}, \dots ,
 s^{t_k - r_k},  \quad {\rm{and}}  \quad   s^{ r_1 } , \dots ,
 s^{ r_k}.
\eeq 
\item[\rm{(vi)}] The subspace  $  W(r, U) $ is hyperinvariant.
\end{itemize}
\end{theorem}

\noindent Proof.
 (i) $\LR$ (ii) It follows from  Lemma \ref{la.bj} that the
two statements are equivalent.
\\
(iv) $\Ra $ (vi) This follows from 
the fact that $ W(r) $ is hyperinvariant. 
\\
(v) $\LR$ (ii) 
Let  $ \tilde {U } \in \mcU $. Then   
 $ W(r,  U) $ and the  quotient space 
$ V / W ( r ,  U  )  $, and also 
  $ W(r,  \tilde {U } ) $ and $ V / W ( r ,  \tilde {U }   )  $,
 have elementary divisors given by \eqref{eq.li2}.
(Note that in
 the  right-hand side of \eqref{eq.wru} 
there may be
 summands of the form
$ \langle u_i  \rangle $ or  $ \langle f ^ {t_i} u_i  \rangle  = 0$.
Thus
\eqref{eq.li2} 
may contain  trivial entries of the form  $ s^0 = 1$.)
\\
(vi) $ \Ra$ (i) Obvious, because of \,$ \rm{Hinv}(V) \subseteq  \rm{Chinv}(V)$.
\\
 (iii) $\Ra$ (iv) 
From 
\,$ \e( u_i  ) = t_i $\, follows 
\[ 
\langle f ^{r_i}  u_i  \rangle   = 
\langle   u_i  \rangle [ f^{ t_i - r_i} ]  \,  \subseteq  \, f ^{r_i} V \cap
V[   f ^{ t _i - r_i} ] .
\]
 Hence
 \,$W (r, U) \subseteq    W (r )$.
We have to show that
the conditions  \eqref{eq.r1}  and \eqref{eq.r2} 
imply the converse inclusion  
\[
    W(  r) =  
f^{r_1}V \cap V[ f^{t_1 - r_1} ]
 \, + \cdots
+ \,
  f^{r_k}V \cap V[ f^{t_k - r_k} ] 
 \subseteq  W(r, U ).
\]
With regard to the decomposition
\,$
    V = \langle  u_1  \rangle  \oplus \cdots \oplus
    \langle  u_k   \rangle $\,
we define 
\[
 D(\mu , \nu ) =  f^{ r_{\nu } }   \langle  u_{\mu}  \rangle
 \,  \cap \, 
\langle  u    _{\mu }   \rangle\,  [f^{ t_{\nu}  - r_{ \nu } } ].
\]
The subspaces 
\,
$    f^ { r_ {\nu} }V  \, \cap  \,  V[ f^ { t_{\nu}  - r_{\nu } } ] $\,
are hyperinvariant.
Therefore Lemma~\ref{la.hdir}(i)  yields 
\[
 f^ { r_ {\nu} }V \cap  V[ f^ { t_{\nu}  - r_{\nu } } ]
=
 \bigoplus  _{ \mu = 1  } ^k   \bigl(
    f^ { r_ {\nu} }V \cap  V [ f^ { t_{\nu}  - r_{\nu } } ]  \cap 
 \langle  u    _{\mu }  \rangle \bigr) 
=
 \bigoplus _{ \mu = 1  } ^k     D(\mu , \nu ). 
\]
Hence 
\beq \label{eq.smm}
     W(r) =
 \, \sum _{\mu, \nu   =  1} ^k   \,  D(\mu , \nu ).
\eeq
Set
$ q(\mu , \nu ) =  \max \{  r_{\nu} ,\,  t_{\mu}  - (  t_{\nu } - r_{\nu} ) \} $.
We have 
\begin{equation*}
 \langle  u _{\mu }   \rangle [f^{ t_{\nu } - r_{ \nu } } ]
=
\begin{cases} 
 \langle  u   _{\mu }   \rangle ,
                                             & 
\text{if \,$  t_{\nu } -  r_{\nu} \ge  t_{\mu}  $},
\\
 f^ { t_{\mu}  - (  t_{\nu } - r_{\nu} ) } \langle  u   _{\mu }  \rangle ,
&
\text{if \,$  t_{\nu } -  r_{\nu} \le  t_{\mu}  $}.
\end{cases}
\end{equation*}
Hence 
\[
  D(\mu , \nu ) 
=  f^ { q(\mu , \nu )   }  \langle  u     _{\mu } \rangle .
\]
Let us show that   $ r_{ \mu} \le   q(\mu , \nu)$ for all $ \mu $.
If
 $ \mu  \geq \nu  $,
 then  \eqref{eq.r2}  implies  
\[
 q(\mu , \nu) = ( t_{\mu}  -   t_{\nu } ) + r_{\nu}  
=  ( t_{\mu}  -  r_{\mu} ) - (   t_{\nu } -  r_{\nu} ) +  r_{\mu}  \ge  r_{\mu} .
\]
If  $ \mu  \leq \nu  $ then   $  t_{\mu}  -   t_{\nu } \le 0 $, and 
therefore  $q(\mu , \nu) =  r_{\nu}$. Hence  \eqref{eq.r1}
 implies   $q(\mu , \nu) \ge  r_{\mu}. $
It follows that 
\[
  D (\mu , \nu) =  f ^{  q(\mu , \nu) } \langle  u _{\mu }  \rangle 
\subseteq 
 f ^{  r_{\nu   } } \langle  u_{\mu }   \rangle 
\subseteq   W(r, U) .
\]
for all $ \mu , \nu $.
 Thus  \eqref{eq.smm} yields 
  $ W(  r)  \subseteq  W( r, U)$.
\\
 (ii) $\Ra $ (iii)
We  modify
the entries of $ U $ 
and 
 replace $ u_k $ by \,$ \tilde{u} _k =  u_{k-1} + u_k  $.
Then 
\,$ \tilde{U} = (u_1 , \dots , u_{k-1} ,  \tilde{u} _k ) \in \mcU$.
Set $ Y_k =  \oplus _{i = 1 }^{k-1}  \langle  f ^{r_i}  u_i   \rangle $.
Then 
$ W (r , U ) =  W (r , \tilde{ U}  )$
implies
\[
    Y_k  \oplus   \langle  f ^{r_k}  u_k  \rangle
=
 Y_k  \oplus   \langle  f ^{r_k} ( u_{k-1} +  u_k  )   \rangle.
\]
From  
\[
   f^{ r_k } u_{k-1 } + f^{r_ k }u_k  \,  \,
\in   \, \,
 \langle  f ^{ r_{k -1} }  u_{k  -1}   \rangle
 \oplus
                             \langle  f ^{r _k}   u_k \rangle
\]
follows 
 \,$  r_{k -1} \le  r_k $.
Proceeding in this manner we obtain
the chain of inequalities in (\ref{eq.r1}).
In order to prove (\ref{eq.r2}) we  start with the  entry of
$ u_1 $ of $ \mcU$ and replace it by
 \,$  u_{1} + f^{ t_2 - t_{1} } u_2 $.
Because of 
\,$\e( u_{1} + f^{ t_2 - t_{1} } u_2 ) =  \e( u_{1} )$
we have 
\,$ \hat{ U } = ( u_1 +  f^{ t_2 - t_{1} } u_2 ,  u_2 , \dots , u_k )
\in \mcU $. 
Set $ Y_1  =  \oplus _{i = 2} ^k  
\langle  f ^{r_i}  u_i  \rangle $.
Then 
$ W (r , U ) =  W (r , \hat{ U}  )$
implies
\[ 
 \langle   f^{ r _1 }  u_1  \rangle \oplus Y_1 
  =
  \langle   f^{ r _1 }   ( u_1 +  f^{ t_2 - t_1 } u_2 )
 \rangle \oplus Y_1 .
\]
From  
\[
  f^{ r _1 }    u_1 +  f^{ r_1 + (t_2 - t_1)  } u_2  
\, \in \,
 \langle   f^{ r _1 }  u_1  \rangle \oplus 
 \langle   f^{ r _2 }  u_2  \rangle 
  \]
follows
\,$  r _2 \le  r_1 + (t_2 - t_1) $, i.e. $ t_1 - r_1 \le t_2 - r_2$,
 such that we  the end up with
\eqref{eq.r2}.
$\sq $

\bigskip 
Let $[ k ]$ denote the greatest integer less than or equal to $k$.
If $ c \in \R $ and 
\,$ 0 < c < 1 $, then
\,$ r = \left( [  c \, t_1  ], \dots ,  [ c \, t_m  ]
\right) $
is  admissible,  and
it is not difficult to verify that $r$
satisfies  \eqref{eq.r1} and \eqref{eq.r2}.
We remark that admissible tuples of the form
$\hat{r}  = ( [\tfrac{1}{2}t_1  ], \dots ,  [\tfrac{1}{2}t_ k ] )$
play a role in the  study   of maximal invariant neutral subspaces
\cite{Ro}.
It follows from Theorem~\ref{thm.qv} that the construction of
such subspaces  is independent of the choice of the
 underlying Jordan basis.

\bigskip

\begin{theorem} \label{thm.int} 
{\rm{(i)}} 
We have 
\beq \label{eq.hdi}
{\rm{Hinv}}(V) =  {\rm{Chinv}}(V) \cap  {\rm{Mark}}(V) . 
\eeq
{\rm{(ii)}}  {\rm{\cite{FHL}}} 
A subspace $ W $ of $V$ is hyperinvariant 
if and only if 
$ W = W(r) $ for some $r $  satisfying  \eqref{eq.r1} and  \eqref{eq.r2}.
\end{theorem} 

\noindent Proof. 
(i) From  Theorem \ref{thm.qv} follows 
\,$  {\rm{Mark}}(V) \cap  {\rm{Chinv}}(V)  \subseteq   {\rm{Hinv}}(V) $.
 The reverse inclusion is 
 \eqref{eq.inch} in  Lemma~\ref{thm.gwa}.
This yields \eqref{eq.hdi}.
Hence  
a  subspace  is hyperinvariant if and only if
it is both characteristic and marked. 
\\
(ii) If  $ W $ is hyperinvariant then $ W $ is marked,
that is $ W = W(r, U )$. Therefore we can apply  
  Theorem \ref{thm.qv}(iv).
It was noted earlier 
 that  $ W(r) \in {\rm{Hinv}}(V)$.
\phantom{00} $\sq $

\bigskip 
We note that
hyperinvariant subspaces can be   characterized completely
by the  distributive law  in Lemma~\ref{la.hdir}.

\begin{theorem} \label{thm.uum}
A subspace $ X  \in {\rm{Inv}}(V)$ is hyperinvariant 
if and only if 
 $ X $ satisfies 
\beq \label{eq.drsss}
  X = (  X \cap V_1 ) \oplus \cdots \oplus  (  X \cap V_q )
\eeq
when
\beq \label{eq.disss}
  V = V_1 \oplus \cdots \oplus V_q, \,\, 
 V_i \in  {\rm{Inv}}(V), \, i = 1, \dots , q.
 \eeq  
\end{theorem}

 \noindent Proof.
Because of  Lemma~\ref{la.hdir} it remains to
prove sufficiency.
Let  $ U = (u_1, \dots , u_k) \in \mcU $ and 
$ \tilde{U}  = ( \tilde{ u } _1, \dots , \tilde{ u } _k) \in \mcU $.
Then  
  \beq \label{eq.tdd}
  V = \langle  u_1  \rangle  \oplus \cdots \oplus
    \langle  u_k   \rangle  =  \langle  \tilde{ u } _1
  \rangle  \oplus \cdots \oplus
    \langle \tilde{ u } _k    \rangle .
 \eeq 
Define 
  \, $ X _i =   \langle  u_i   \rangle  \cap X $ and  
 $ \tilde{X} _i =
    \langle  \tilde{ u }_i   \rangle  \cap X $,
$i = 1, \dots , k$.
Then 
\,  $ X _i
 =  \langle  f ^{r_i}  u_i   \rangle $ and 
 $ \tilde{X} _i =
 \langle  f ^{ \tilde{r}_i} \tilde{ u}_i     \rangle $
for some $ r_i,  \tilde{r}_i $.
Set \,$ r = ( r_1, \dots , r_k ) $\, and \,$ \tilde{r} =
( \tilde{r}_1 , \dots ,  \tilde{r}_k ) $.
In \eqref{eq.tdd} we have two
direct sums of the form \eqref{eq.disss}.
Hence the assumption 
\eqref{eq.drsss}  implies \,$X = W(r,U) = W ( \tilde{r} , \tilde{U} )   $.
 We can pass
 from  $ U $ to  $ \tilde{U} $  in at most $k$ steps,
changing a single entry at each step. 
Suppose we replace $ u_k$ in $U $ by  $  \tilde{u}_k $.
Then 
\,$ \hat{ U}   = (u_1,   \dots , u_{k- 1} ,
  \tilde{u}_k ) 
\in \mcU $,
and \,$
  V = \langle  u_1  \rangle  \oplus \cdots
\langle  u_{k-1}  \rangle  \oplus \langle  \tilde{u}_k   \rangle $.
Set $ Y_k = \oplus _{i = 1 } ^{k-1}  
\langle f^{r_i}   u_i   \rangle  $.
Then  
 \[  
  X = Y_k \oplus  
     \langle  f ^ { \tilde{r} _k}  \tilde{u}_k   \rangle 
= 
 Y_k  \oplus        \langle  f ^ { r _k}  u_k   \rangle .   
 \]
  Considering the elementary divisors of $ V/ X $ we
 deduce $ \tilde{r}_k = r_k $, 
 and at the end we obtain  \,$ r =  \tilde{r}$, and therefore
$  W (r, U) =  W (r, \tilde{U})$.
We conclude that $ X = W (r, U) $ is 
independent of the choice of the generator tuple  $U $. Hence 
$ X $ is hyperinvariant. 
$ \sq $

Let us  reexamine 
 Example  \ref{ex.fu22} and consider a field $K$ of charac\-teristic 
different from $2 $.

\noindent {\bf{Example  \ref{ex.fu22} (continued).}}  
Let  ${\rm{char}} K \ne 2 $. Then 
$ \gamma : (e_1, e_2) \mapsto ( 2 e_1, e_2)   $ determines
an $f$-automorphism. For     $ Z = \langle  e_1 + e_3 \rangle  $
we have 
$ \gamma (Z ) =  \langle 2 e_1 + e_3 \rangle  \ne Z $.
Hence in this case 
 \,$ Z \in {\rm{Inv}}(V)$\,  
is not characteristic. 

To identify the
 characteristic subspaces we screen $  {\rm{Inv}}(V) $. 
Note that 
\begin{multline*} 
 {\rm{Aut}}_f(V) =
\big\{  
\alpha :  (e_1, e_2 ) \mapsto (  a e_1 + b e_ 4,
 c e_2 + d e_3 + g e_4 + h e_1 ),
 \\  a, b , c , d ,g ,h \in K, 
 a \ne  0, c \ne  0 
 \big\}.
\end{multline*}
The nonzero cyclic 
 subspaces
 are of the form 
$  \langle e_2  + c e_1   \rangle $,
 $  \langle e_3  + c e_1   \rangle $,
and  $ \langle a  e_4  + c e_1   \rangle $, $a, c \in K$, 
 $ (a, c ) \ne (0,0) $.
Only  $\langle  e_3    \rangle =   f V $ and 
 $
 \langle e_4  \rangle =  f^2 V $ are  characteristic.
Moreover,  $X$ is  a direct sum of two cyclic subspaces
if and only if $ X \in \{  V , \langle  e_3    \rangle   
\oplus  \langle e_1  \rangle  =  V [f^2  ],   \langle  e_4    \rangle   
\oplus  \langle e_1  \rangle  =  V [f ] \} $. 
These  three  subspaces are  characteristic. 
We find \,$ {\rm{Hinv}}(V) = \{ 0, \,  f V , \, f^2 V , \,
 V[f ] , \,
V[f ^2 ] , \, V\} $.
Hence 
$ {\rm{Hinv}}(V) =  {\rm{Chinv}}(V) $.
The example  is a special case of  the following    
 general result (see also \cite[p.\,67]{Kap}).

\begin{theorem} \label{thm.k2}  
If \,  
  $ | K | > 2 $\,  
then each characteristic subspace of   $V$   is hyperinvariant,
i.e. 
$ 
  {\rm{Chinv}}(V)  =   {\rm{Hinv}}(V) $.
\end{theorem}

\noindent Proof. 
Because of Lemma~\ref{la.ds} it suffices to consider 
the case where $ f $ has only one eigenvalue.
We can assume $ f ^n = 0$. 
If $ | K | > 2 $  and
 $ X $ is characteristic  then 
 it follows from 
  Lemma \ref{la.hdir}(ii)  
that \eqref{eq.dis} implies \eqref{eq.drs}. 
Therefore, according to  Theorem \ref{thm.uum}, 
the subspace $X $ is hyperinvariant. 
$\sq$

\medskip 
In the case of  vector spaces   over  $K = \Z _2 $ 
it is an open problem to describe 
all subspaces  that are characteristic
without being hyperinvariant.

\medskip 
\noindent 
{\bf{
 Acknowledgment:}} We are  indebted to L. Rodman
for a valuable remark.

\end{document}